%% file: VarCat18.tex
\documentclass{amsart}
\usepackage{graphicx}
\usepackage{verbatim}
\newtheorem{thm}{Theorem}[section]
\newtheorem{cor}[thm]{Corollary}
\newtheorem{lem}[thm]{Lemma}
\newtheorem{prop}[thm]{Proposition}
\newtheorem{conj}[thm]{Conjecture}
\theoremstyle{definition}

\theoremstyle{remark}
\newtheorem{rem}[thm]{Remark}

\numberwithin{equation}{section}

\newcommand{\set}[1]{\left\{#1\right\}}
\newcommand{\Real}{\mathbb R}

\newcommand{\func}[1]{\ensuremath{\mathrm{#1} \:} }

\newcommand{\Div}[0]{\func{div}}

\newcommand{\re}[0]{\func{Re}}
\newcommand{\im}[0]{\func{Im}}

\newcommand{\xX}[0]{\mathbf{x}}
\newcommand{\eE}[0]{\mathbf{e}}
\newcommand{\FF}[0]{\mathbf{F}}

\title{A Variational Characterization of the Catenoid}
\subjclass[2000]{53A10}
\author{Jacob Bernstein and Christine Breiner}
\address{Dept. of Math,
Stanford University, Stanford, CA 94305, USA}

\email{jbern@math.stanford.edu}
\address{Dept. of Math, Massachusetts Institute of
Technology, Cambridge, MA  02139, USA} \email{breiner@math.mit.edu}
\thanks{The authors were supported respectively by the NSF grants DMS-0902721 and DMS-0902718.}

\begin{document}
\begin{abstract}
In this note, we use a result of Osserman and Schiffer \cite{OS} to give a variational characterization of the catenoid.  Namely, we show that subsets of the catenoid minimize area within a geometrically natural class of minimal annuli.  To the best of our knowledge, this fact has gone unremarked upon in the literature.  As an application of the techniques, we give a sharp condition on the lengths of a pair of connected, simple closed curves $\sigma_1$ and $\sigma_2$ lying in parallel planes that precludes the existence of a connected  minimal surface $\Sigma$ with $\partial \Sigma=\sigma_1\cup\sigma_2$.
\end{abstract}
\maketitle
\section{Introduction}
Recall that the catenoid is the minimal surface of revolution surface given by
\begin{equation*}
Cat=\set{x_1^2+x_2^2=\cosh^2 x_3}\subset \Real^3.
\end{equation*}
The catenoid was discovered by Euler in 1744 and is one of the classic examples in the theory of minimal surfaces and, more broadly, in the calculus of variations.  For instance, a sequence of homothetic blow-downs of $Cat$ provides the simplest model of the failure of smooth convergence for a sequence of minimal surfaces.
  Due to the invariance of the minimal surface equation under rigid motions and homotheties of $\Real^3$, $Cat$ sits within a six dimensional family of catenoids which we henceforth denote by $\mathcal{C}$.  In other words, $C\in \mathcal{C}$ if $C$ can be obtained from $Cat$ by a rigid motion composed with a homothety.

The catenoid (or rather $\mathcal{C}$) has been characterized in many ways.  We list some notable results:  In the spirit of Euler, O. Bonnet, in the mid-nineteenth century, showed that the catenoid and plane are the only minimal surfaces of revolution.  More recently, in \cite{SCH}, R. Schoen showed that the catenoid is the unique complete, embedded, minimal surface with finite total curvature and two ends.  In a similar vein, but by very different techniques, F. J. L\'opez and A. Ros showed that the catenoid and plane are the unique complete, embedded minimal surfaces of finite total curvature and genus-zero \cite{LopezRos}.  Additionally, building on work of D. Fischer-Colbrie \cite{FC}, L\'opez and Ros characterized the catenoid as the unique complete, embedded minimal surface in $\Real^3$ of Morse index one -- see \cite{LopezRos2}.    Finally, we note that in \cite{Pyo} pieces of the catenoid are shown to be the only minimal annuli in a slab that meet the boundary of the
  slab in a constant contact angle -- a fact that will be relevent in this note.  It bears mentioning that work of P. Collin \cite{Co} and T.H. Colding and W. P. Minicozzi \cite{CY} allows one to replace the geometric assumption of finite total curvature in \cite{LopezRos} and \cite{SCH} with the weaker topological assumption that the surfaces are of finite topology -- that is, diffeomorphic to finitely punctured compact surfaces.

 In this note, we characterize the catenoid as the unique minimal surface which minimizes area within a geometrically natural class of minimal annuli.  This turns out to be a simple consequence of the proof due to R. Osserman and M. Schiffer \cite{OS} of the isoperimetric inequality for minimal annuli in $\Real^3$.

Let us now describe in what sense the catenoid minimizes area.  Fix two parallel planes $P_-,P_+\subset \Real^3$ with $P_-\neq P_+$ and denote by $\Omega$ the open slab between them.    We remark that for any plane $P\subset \Omega$, $P$ must be parallel to $P_-$. 
Let us denote by $\mathcal{M}(\Omega)$ the class of smooth minimal surfaces spanning $P_-$ and $P_+$.  That is, $\Sigma\in \mathcal{M}(\Omega)$ if $\Sigma\subset \Omega$ may be parameterized by a conformal, harmonic immersion $\mathbf{F}: M\to \Real^3$ so that $b\Sigma:=\overline{\Sigma}\backslash \Sigma \subset \partial \Omega=P_+\cup P_-$.    Here $M$ is an open orientable surface.  Notice that an element $\Sigma \in \mathcal{M}(\Omega)$ may have arbitrarily bad behavior as one approaches $\partial \Omega$; however, if $\Sigma$ has the structure of a surface with boundary then $\partial \Sigma=b\Sigma$.
The class $\mathcal{M}(\Omega)$ is too broad for our methods and so we will restrict attention to the subclass $\mathcal{A}(\Omega)\subset \mathcal{M}(\Omega)$ consisting of embedded annuli.  Precisely, $\Sigma \in \mathcal{A}(\Omega)$ if, in addition to lying in $\mathcal{M}(\Omega)$, $\Sigma$ may be parameterized by an annulus; i.e.  there is a smooth embedding $\mathbf{F}: (0,1)\times \mathbb{S} \to \Real^3$ with image $\Sigma$.  It bears mentioning that in \cite{MeeksWhite}, by using global analysis techniques, W. Meeks and B. White have shown that the subset of $\mathcal{A}(\Omega)$ consisting of surfaces with $b\Sigma$ a pair of $C^{2,\alpha}$ (planar) convex curves has the structure of a contractible Banach manifold.

Recall that given a pair of connected, simple closed curves $\sigma_+\in P_+$ and $\sigma_- \in P_-$ there need not exist $\Sigma \in \mathcal{M}(\Omega)$ with $\partial \Sigma=\sigma_+ \cup \sigma_-$.  Indeed, for $\sigma_\pm$ of sufficiently small length relative to the size of the slab, the existence of such a surface would violate the isoperimetric inequality for minimal surfaces. In Theorem 4 of \cite{OS}, Osserman and Schiffer give a sharp condition on the lengths of the $\sigma_\pm$ that precludes the existence of a minimal annulus spanning the $\sigma_\pm$.  In Section \ref{NoSpanSec}, we refine their result and show that their condition actually precludes the existence of \emph{any} connected minimal surface spanning the curves and prove, in addition, an interesting rigidity result.

It is important to emphasize that our class of surfaces consists of \emph{minimal} surfaces.  If one were to consider classes of arbitrary surfaces spanning $P_-$ and $ P_+$ then the infimum of area would be zero -- as can be seen by considering thin cylinders. Similarly, the surfaces in $\mathcal{M}(\Omega)$ are, in general, not stationary for area with respect to variations moving their boundary.  In particular, we are not considering the free boundary problem as usually formulated for minimal surfaces.

We claim that for some $C\in \mathcal{C}$, $C\cap \Omega\in \mathcal{A}(\Omega)$ has less area than any other surface in $\mathcal{A}(\Omega)$.   More precisely, denote by:
\begin{equation*}
Cat_{MS}=\set{x_1^2+x_2^2=\cosh^2 x_3}\cap \set{1-x_3\tanh x_3> 0}\subset Cat
\end{equation*}
 the maximally symmetric marginally stable piece of the catenoid.  That is $Cat_{MS}$ is stable but any domain in $Cat$ strictly containing $Cat_{MS}$ is unstable.  Recall a minimal surface  is \emph{stable} if no compactly supported infinitesimal deformation decreases area; it is \emph{unstable} if there is a compactly supported infinitesimal deformation decreasing area.  Marginally stable surfaces are on the boundary between these two classes.  We show that $Cat_{MS}$ provides the model for least area surfaces in $\mathcal{A}(\Omega)$.  
\begin{thm}\label{MainThm}
Let $P_-$ and $P_+$ be distinct parallel planes in $\Real^3$ and let  $\Omega$ be the open slab between them. 
Let $C_{MS}\in \mathcal{A}(\Omega)$ be the (unique up to translations parallel to $P_\pm$)  minimal surface in $\mathcal{A}(\Omega)$ obtained from rigid motions and homotheties of $Cat_{MS}$.  
Then for any $\Sigma \in \mathcal{A}(\Omega)$:
\begin{equation}
\mathcal{H}^2(\Sigma)\geq \mathcal{H}^2(C_{MS})
\end{equation}
 with equality if and only if $\Sigma$ is a translate of $C_{MS}$.  
\end{thm}
Here $\mathcal{H}^k$ denotes $k$-dimensional Hausdorff measure.  We point out that the two disks $D_\pm \subset P_\pm$ with $\partial D_+ \cup \partial D_- = \partial C_{MS}$ satisfy $\mathcal{H}^2(D_-\cup D_+)< \mathcal{H}^2(C_{MS})$.
That is $C_{MS}$ is not an area minimizer with respect to its boundary even though it does minimize area in the class of spanning minimal annuli $\mathcal{A}(\Omega)$. 

The restriction to $\mathcal{A}(\Omega)$ in Theorem \ref{MainThm} rather than $\mathcal{M}(\Omega)$ is necessitated by our argument.  However, it seems reasonable to believe that an area minimizer, $\Sigma_0 \in\mathcal{M}(\Omega)$ should have $b\Sigma_0$ rather nice -- for instance consisting of convex planar curves. 
 Hence, in light of the embeddedness results of T. Ekholm, B. White and D. Wienholtz \cite{EBW} and a long standing conjecture of W. Meeks on the non-existence of postive genus surfaces in $\mathcal{M}(\Omega)$ bounded by convex curves (see Conjecture 3.10 of \cite{MeeksConj}), it is natural to conjecture:
\begin{conj}
Let $P_-$ and $P_+$ be distinct parallel planes in $\Real^3$ and let  $\Omega$ be the open slab between them. 
Let $C_{MS}\in \mathcal{M}(\Omega)$ be the (unique up to translations parallel to $P_\pm$)  minimal surface in $\mathcal{M}(\Omega)$ obtained from rigid motions and homotheties of $Cat_{MS}$.  
Then for any $\Sigma \in \mathcal{M}(\Omega)$:
\begin{equation}
\mathcal{H}^2(\Sigma)\geq \mathcal{H}^2(C_{MS})
\end{equation}
 with equality if and only if $\Sigma$ is a translate of $C_{MS}$.  
\end{conj}

 Returning to the more restricted setting  of $\mathcal{A}(\Omega)$, we note that Theorem \ref{MainThm} is a simple consequence of a more general area minimization property of the catenoid.  Indeed, minimal annuli have a natural scale which may be computed as the length of the flux vector associated to the generator of the homology group -- we refer to Section \ref{FluxSec} for precise definitions.  Normalizing with respect to this scale gives an area lower bound:
\begin{thm}\label{SecondThm}
Let $P_-=\set{x_3=h_-}$ and $P_+=\set{x_3=h_+}$ be distinct parallel planes in $\Real^3$ with $h_-<h_+$ and let $\Omega$ be the open slab between them.   Fix $\Sigma \in \mathcal{A}(\Omega)$.  
Let $P_0=\set{x_3=h_0}\subset \bar{\Omega}$ denote the plane that satisfies:
\begin{equation*}
\mathcal{H}^1 (\Sigma \cap P_0)=\inf_{t\in (h_-,h_+)} \mathcal{H}^1(\Sigma_t).
\end{equation*}
Here $\Sigma_t=\Sigma\cap \set{x_3=t}$ and $\mathcal{H}^1 (\Sigma \cap P_+)$ is defined as $\liminf_{t\nearrow h_+ } \mathcal{H}^1 (\Sigma_t)$ and likewise for $\mathcal{H}^1 (\Sigma \cap P_-)$.
Let $F_3$ denote the vertical component of $Flux(\Sigma)$. If $C$ is the vertical catenoid with $Flux(C)=(0,0,F_3)$ and symmetric with respect to reflection through the plane $P_0$ then:
\begin{equation}
\mathcal{H}^2(\Sigma)\geq \mathcal{H}^2(C\cap  \Omega)
\end{equation}
 with equality if and only if $\Sigma$ is a translate of  $C\cap \Omega$.  
\end{thm}

The proof of Theorem \ref{SecondThm} will comprise the bulk of this note.  Let us first use it to prove Theorem \ref{MainThm}:
\begin{proof}
By Theorem \ref{SecondThm}, the least area surface must be a piece of a catenoid.  Up to a rescaling and rigid motion we may take $\Omega=\set{-1<x_3<1}$ and so may restrict attention to subsets of vertical catenoids.   The space of these catenoids are parameterized by $\lambda>0$ and $t$ where
\begin{equation*}
Cat_{\lambda,t}=\lambda Cat+t\mathbf{e}_3.
\end{equation*}
Set
\begin{equation*}
A(\lambda,t)=\mathcal{H}^2(Cat_{\lambda,t}\cap \Omega).
\end{equation*}
We will check that this value is minimized at $\lambda_0, t_0$ chosen so $Cat_{\lambda_0, t_0} \cap \Omega=\lambda_0 Cat_{MS}$.

Let us parameterize a scale $\lambda$ vertical catenoid by 
\begin{equation*}
\mathbf{F}_\lambda (h, \theta)=(\lambda \cosh \frac{h}{\lambda}\cos \theta , \lambda\cosh \frac{h}{\lambda} \sin \theta, h).
\end{equation*}
Then we have by the area formula:
\begin{align*}
A( \lambda,t)&=\lambda \int_{-1-t}^{1-t} \int_0^{2\pi}  \cosh^2 \frac{h}{\lambda} d\theta dh\\
   					&=2\pi \lambda^2 \left.\left( \frac{1}{4} \sinh \frac{2h}{\lambda}+\frac{h}             {2\lambda}\right)\right|_{h=-1-t}^{h=1-t} \\
   					&=\frac{\lambda^2 \pi}{2} \left(\sinh \frac{2(1-t)}{\lambda}-\sinh  \frac{2(-1-t)}{\lambda} \right) +2\pi\lambda\\
   					&={\lambda^2 \pi}\cosh \frac{2t}{\lambda}  \sinh \frac{2}{\lambda}+2\pi\lambda.
 \end{align*}
It is elementary to check that $A(\lambda,t)\to \infty $ as $|(\ln \lambda,t)|\to \infty$, so in order to find the least area catenoid we look for critical points of $A$.
 We expect these to occur for catenoids that are symmetric about the plane $\set{x_3=0}$. Indeed,
$\partial_t A=2\pi \lambda \sinh\frac{2t}{\lambda} \sinh \frac{2}{\lambda}$ and this equals zero if and only if $t=0$ and hence $t_0=0$.  Thus, we need only find $\lambda$ so
\begin{equation*}
0=\partial_\lambda A(\lambda,0)=2\pi \lambda \sinh \frac{2}{\lambda}- 2\pi\cosh \frac{2}{\lambda}+2\pi
\end{equation*}
Using $\cosh x=\sqrt{1+\sinh^2 x}$,  this equation is equivalent to solving:
\begin{equation*}
(\lambda^2-1) \sinh^2 \frac{2}{\lambda} +2\lambda \sinh \frac{2}{\lambda} =0
\end{equation*}
which has a unique solution $\lambda=\lambda_0\approx 0.833$ determined by
\begin{equation}\label{lambdanot}
\frac{2\lambda}{1-\lambda^2} = \sinh \frac{2}{\lambda}.
\end{equation}
Uniqueness follows from properties of the two functions in \eqref{lambdanot} on the domain $\lambda \in [0, 1)\cup (1, \infty)$.  First, observe that the right hand side is always positive on this domain while the left hand function is only positive on $[0,1)$.  Second, $\sinh \frac{2}{\lambda}$ is strictly decreasing on $[0,1)$ while $\frac{2\lambda}{1-\lambda^2}$ is strictly increasing. 
Finally, $\lim_{\lambda\searrow 0}  \sinh \frac{2}{\lambda}-\frac{2\lambda}{1-\lambda^2}=\infty$ while  $\lim_{\lambda\nearrow 1}  \sinh \frac{2}{\lambda}-\frac{2\lambda}{1-\lambda^2}=-\infty$. 

Lastly, we must verify that $Cat_{\lambda_0, t_0}\cap \Omega=\lambda_0 Cat_{MS}$.  To that end we note that
standard properties of hyperbolic functions give that
\begin{equation*}
\sinh \frac{1}{\lambda_0}=\frac{\lambda_0}{\sqrt{1-\lambda_0^2}}
\end{equation*}
and so
\begin{equation*}
\tanh \frac{1}{\lambda_0}=\lambda_0.
\end{equation*}
Hence, $1-\frac{1}{\lambda_0} \tanh \frac{1}{\lambda_0}=0$.  But this implies that $1-\frac{1}{\lambda_0} \tanh \frac{x_3}{\lambda_0}>0$ on $ Cat_{\lambda_0, t_0}\cap \Omega$.  That is, $\frac{1}{\lambda_0} (Cat_{\lambda_0, t_0}\cap \Omega)= Cat_{MS}$.
\end{proof}

\subsection*{Acknowledgments}
We would like to thank David Hoffman for his interest and his many helpful comments.  We
also thank Brian White for suggesting a simplification of the proof of Lemma \ref{AltLem}.
\section{Convexity of the length of level sets}
Before we proceed with the proof of Theorem \ref{SecondThm} we must recall some important definitions.  We also state the result of Osserman and Schiffer regarding the convexity of the length of certain families of curves in minimal annuli in $\Real^3$. 
\subsection{The Flux Vector} \label{FluxSec}
  For the purposes of this discussion we assume that $\Sigma\in \mathcal{A}(\Omega)$ for some $\Omega$.  Fix an orientation of $\Sigma$ and let $\gamma\subset \Sigma$ be a simple $C^1$ closed curve in $\Sigma$ on which we also fix an orientation.  
Our choices of orientation give rise to a normal vector field in $\Sigma$ along $\gamma$ which we denote by $\nu$.  We always think of the vectors $\nu$ as vectors in $\Real^3$. The \emph{flux} of $\gamma$ is defined to be the vector: 
\begin{equation}
Flux(\gamma)=\int_{\gamma} \nu  \; d\mathcal{H}^1\in \Real^3.
\end{equation}
 As $\nu\cdot \mathbf{e}_i =\nu\cdot \nabla_\Sigma x_i$ and on a minimal surface $\Delta_\Sigma x_i=0$, the divergence theorem implies that the flux of a curve depends only on its homology class. In particular, for a minimal annulus, $\Sigma$, we may associate a vector $Flux(\Sigma)$, by choosing $\gamma$ so that $[\gamma]$ is  a generator of $H_1(\Sigma)$ and setting $Flux(\Sigma)=Flux(\gamma)$;  up to a reflection through the origin, $Flux(\Sigma)$ is independent of the choice of orientation of $\Sigma$ and of $\gamma$.  In the sequel, we will consider 
\begin{equation*}
F_3(\Sigma) = Flux(\Sigma) \cdot \mathbf{e}_3,
\end{equation*}
the vertical component of the flux of the minimal annulus $\Sigma$.  We always choose orientations so that $F_3 \geq 0$.

An important property of the flux is that it sets a natural scale for a minimal annulus.  Namely, suppose that $\Sigma$ is a minimal annulus and $\Sigma'=\lambda \Sigma$ is the annulus obtained by homothetically scaling $\Sigma$ by $\lambda>0$. Then one computes $Flux(\Sigma')=\lambda Flux(\Sigma)$.    In particular, the flux allows one to distinguish between catenoids of differing scales.
A more subtle property of the flux is that it also helps to set a natural conformal scale for  elements of $\mathcal{A}(\Omega)$.  
More precisely consider $\Sigma\in \mathcal{A}(\Omega)$.  By the uniformization theorem there is a conformal diffeomorphism $\psi$ between $\Sigma$ and a flat open cylinder, $(h_-, h_+)\times \mu \mathbb{S}$ where here $(h_-, h_+)$ denotes a (possibly infinite) interval and $\mu \mathbb{S}$ denotes the circle of radius $\mu$.   Moreover, the ratio between $|h_+-h_-|$ and $\mu$ is determined by $\Sigma$.  We claim this ratio is actually determined only by $\Omega$ and $Flux(\Sigma)$.  

In order to show this we first need the following fact:
\begin{lem}\label{TransverseLem}
Let $\Sigma\in \mathcal{A}(\Omega)$ and suppose that $\bar{\Sigma}$ has the structure of a surface with boundary and that $\partial\Sigma$ is smooth.  Then for any plane $P\subset \bar{\Omega}$, $P$ meets $\Sigma$ transversally.
\end{lem}
\begin{proof}
Using an ambient rigid motion and homothety, we take $\Omega=\set{-1<x_3<1}$.
As $\partial \Sigma$ is smooth, standard boundary regularity results imply that $\Sigma$ may be viewed as a smooth surface with boundary. 
We first show that $\Sigma$ meets the planes $P_1=\set{x_3=1}$ and $ P_2=\set{x_3=-1}$ transversally.  To that end we note that as 
$\partial \Sigma$ is smooth there is a uniform constant $r_0>0$ so that for any point $p\in \partial \Sigma$, the inner and outer 
osculating circles have radius greater than $r_0$.  Moreover, there is a uniform bound on the ratio between intrinsic and extrinsic 
distance between points of $\partial \Sigma$. Hence, there exists $r_1$ with $0<r_1<r_0$ so for any $p \in P_1\cap \partial \Sigma$ the following holds:  
there are circles in $P_1$ denoted by $C_{in}(p)$ and $C_{out}(p)$, both of radius $r_1$, such that $C_{in}(p)$ lies within 
$\partial \Sigma\cap P_1$ (thought of as a plane curve in $P_1$) while $C_{out}(p)$ lies outside $\partial \Sigma \cap P_1$.  Moreover, both circles $C_{in}(p), C_{out}(p)$
 meet $\partial \Sigma$ only at $p$. A similar result holds for 
$p\in \partial \Sigma \cap P_2$.  Without loss of generality we consider only $p\in \partial \Sigma\cap P_1$.

Now let $Cat^+=Cat \cap \set{x_3\geq 0}$.  Denote by $Cat_{in}(p)$ the set obtained from $Cat^+$ by 
translations and homotheties so that $\partial 
Cat_{in}(p)=C_{in}(p)$ and let $Cat_{out}(p)$ be defined in an analogous manner; 
notice that both $Cat_{in}(p)$ and $Cat_{out}(p)$ are disjoint from $\Omega$.  
Denote by $Cat_{in}'(p)$ the surfaces obtained by scaling $Cat_{in}(p)$ by $\frac{1}{2}$ about the center of $C_{in}(p)$ and define  
$Cat_{out}'(p)$ similarly.  By the strict maximum principle and the definition 
of $\mathcal{A}(\Omega)$ we have that $\overline{\Sigma} \cap \partial \Omega=\partial \Sigma$.  Hence, there is a $\delta>0$ so that for all $p$ both 
$Cat_{in}'(p)$ and $Cat_{out}'(p)$ can be translated along their axes by $\delta$ into $\Omega$ while 
remaining disjoint from $\Sigma$ (that is translated in the direction $-\textbf{e}_3$).  Let us denote by $Cat_{in}''(p)$ and $Cat_{out}''(p)$ the surfaces resulting from this translation and let $\lambda Cat_{in}''(p)$ denote the result of scaling $Cat_{in}''(p)$ by 
$\lambda>0$ about the center of $\partial Cat_{in}''(p)$ and similarly for $\lambda Cat_{out}''(p)$.  As 
$\lambda \to 0$ both $\lambda Cat_{in}''(p)$ and $\lambda Cat_{out}''(p)$ converge to a plane $P''(p)\subset \Omega$ 
which must meet $\Sigma$ as otherwise $\partial \Sigma \cap P_1=\emptyset$ which is impossible by the convex hull 
property and our definition of $\mathcal{A}(\Omega)$.   Hence, there are $0<\lambda_{in}
(p)<1$ and $0<\lambda_{out}(p)<1$ so that,  for $\lambda<\lambda_{in}(p)$, $\lambda  Cat_{in}''(p)$ meets $\Sigma$ 
but, for $\lambda>\lambda_{in}(p)$, $\lambda  Cat_{in}''(p)$ is disjoint from $\Sigma$; and the same for 
$\lambda_{out}(p)$ with respect to $\lambda Cat_{out}''(p)$.   As a consequence, $\lambda_{in}(p)Cat_{in}''(p)$ is 
disjoint from $\Sigma\backslash \partial \Sigma$ but meets $\partial \Sigma$ and it must do so precisely at $p$ and the same is true for $\lambda_{out}(p)Cat_{out}''(p)$.   By the boundary maximum principle we 
then see that the normal to $\Sigma$ at $p$ cannot be orthogonal to $P_1$ and hence $P_1$ meets $\Sigma$ transversally as do all planes $\set{x_3=t}$ for $1-\epsilon<t\leq 1$. A similar result holds for $P_2$.  Thus, for all planes $P\subset \Omega$ near $P_1$ or $P_2$,  $P\cap \Sigma$ consists of a single smooth simple closed curve.

Now let $f=x_3$ be the function whose level sets are planes in $\Omega$.  As $\Sigma$ is minimal $f$ is a harmonic function on $\Sigma$ and so has no local maxima or minima. In particular, at any critical points of $f$ the vector field $\nabla f$  has negative index.  By our previous discussion $f$ has no critical points near $\partial \Sigma$ and, moreover, $\nabla f$ is transverse to $\partial \Sigma$.  As $\Sigma$ is an annulus, the Hopf index theorem then implies that $f$ has no critical points.  
\end{proof}
As a consequence, $Flux(\Sigma)$ and $\Omega$ determine the cylinder with which $\Sigma$ is conformally equivalent:
\begin{cor} \label{DiffeoCor}
Suppose that $\Omega=\set{h_-<x_3<h_+}$ is an open slab.  For any $\Sigma \in \mathcal{A}(\Omega)$ there is a conformal diffeomorphism
\begin{equation*}
\psi:\Sigma \to (h_-,h_+)\times \mu \mathbb{S}
\end{equation*}
given by $\psi(p)=(x_3(p), x_3^*(p))$.  Moreover, $\mu=\frac{1}{2\pi} F_3(\Sigma )>0$.
\end{cor}
\begin{proof}
Let $J$ be the almost-complex structure on $\Sigma$ arising from the metric and some choice of orientation. As $x_3$ is a harmonic function,  $dx_3$ is a harmonic one form.  Moreover, $dx_3^*=dx_3\circ J$, the conjugate differential,  is also harmonic. In general, $dx_3^*$, while closed, will never be exact.  Indeed, for a general closed curve $\gamma$ in $\Sigma$:
\begin{equation*}
Flux(\gamma)\cdot \mathbf{e}_3=\int_{\gamma} dx_3^*.
\end{equation*}
Thus, integrating $dx_3^*$ gives a map $x_3^*: \Sigma \to \mathbb{R}/F_3\mathbb{Z}$ where $F_3=F_3(\Sigma)\geq 0$.   We will see that we must have $F_3> 0$ and so $ \mathbb{R}/F_3\mathbb{Z}=\mu \mathbb{S}$. 
As $dx_3$ and $dx_3^*$ have the same length and are orthogonal, if we set $\psi=(x_3, x_3^*)$ then $\psi$ is a conformal map. 

By Sard's theorem, for each $\epsilon>0$ there is an $0<\epsilon' <\epsilon$ so that both $\set{x_3=h_-+\epsilon'}$ and $\set{x_3=h_+-\epsilon'}$ meet $\Sigma$ transversally.  Moreover, as $b\Sigma \subset \partial \Omega$, each of the finitely many components of $\set{h_-+\epsilon'\leq x_3\leq h_--\epsilon'}\cap \Sigma$ has the structure of a surface with boundary.  By the convex hull property of minimal surfaces there is exactly one such component $\Sigma_{\epsilon'}$ and it is an annulus.
 Lemma \ref{TransverseLem} implies that for $h_-+\epsilon'<t<h_+-\epsilon'$, each plane $\set{x_3=t}$ meets $\Sigma_{\epsilon'}$, and hence $\Sigma$, transversally.  In particular, $dx_3$ does not vanish on  $\Sigma_{\epsilon'}$  and hence  $\psi$  restricted to $\Sigma_{\epsilon'}$ is a local diffeomorphism.    As each level set of $x_3=t$ for $h_-+\epsilon'\leq t\leq h_+-\epsilon'$ is connected and $dx_3^*$ does not vanish, $\psi$ is injective on these level sets and hence the restriction of $\psi$ to $\Sigma_{\epsilon'}$ is injective.  In addition, it is then clear that $F_3\neq 0$.  Taken together it follows that $\psi$ restricts to a conformal diffeomorphism between $\Sigma_{\epsilon'}$ and $ (h_-+\epsilon',h_+-\epsilon')\times \mu \mathbb{S}$.  As $\epsilon$ may be taken as small as we like, the result is shown.
\end{proof} 

\subsection{Osserman and Schiffer's Result}
We now record the convexity result of Osserman and Schiffer from \cite{OS} that we will use. This result was a key step in their proof -- also in \cite{OS} -- of the sharp isoperimetric inequality for doubly connected minimal surfaces in $\Real^3$.  We point out  that the restriction to $\Real^3$ comes from their use of the Weierstrass representation in order to prove the convexity result.  Roughly speaking, Osserman and Schiffer show that when a minimal annulus $\Sigma\subset\Real^3$ is conformally parametrized by an annulus $A$ in the complex plane, then the length of the images in $\Sigma$ of the circles foliating $A$ satisfy a convexity condition that is sharp on catenoids and planar annuli.  Precisely,
\begin{lem}\label{ShiffOssLem}
Let $A_{r,R}=\set{z: r<|z|<R}\subset \mathbb{C}$ and suppose that $\mathbf{F}: A_{r,R}\to \Real^{3}$ is a conformal harmonic immersion (so in particular the image of $\mathbf{F}$ is a minimal surface).  If we let $\sigma_\rho$ be the image of $|z|=\rho$ under $\mathbf{F}$ and define:
\begin{equation}
L(t)=\mathcal{H}^1(\sigma_{e^t})
\end{equation}
then
\begin{equation}
L''(t)\geq L(t)
\end{equation}
with equality if and only if $\mathbf{F}$ maps into a planar annulus or into a piece of a catenoid bounded by coaxial circles in parallel planes.
\end{lem}
For the sake of completeness we sketch Osserman and Schiffer's proof in Appendix \ref{OSApp}.
Rather than using Lemma \ref{ShiffOssLem} directly we use the following corollary:  
\begin{cor}\label{OSCor}
Suppose that $\Omega=\set{h_-<x_3<h_+}$ and $\Sigma \in \mathcal{A}(\Omega)$.  Set $\Sigma_t=\Sigma\cap \set{x_3=t}$. Then for $t\in (h_-,h_+)$:
\begin{equation*}
\frac{d^2}{dt^2} \mathcal{H}^1(\Sigma_t) \geq \frac{(2\pi)^2}{F_3(\Sigma)^2} \mathcal{H}^1(\Sigma_t)
\end{equation*}
with equality if and only if $\Sigma$ is a piece of a vertical catenoid $C$.  
\end{cor}
\begin{proof}
By Corollary \ref{DiffeoCor}, $F_3\neq 0$ and there is a conformal diffeomorphism
\begin{equation*}
\psi :(h_-,h_+)\times \mu \mathbb{S}\to \Real^3
\end{equation*}
with image $\Sigma$ and so that $\Sigma_t$ is the image of $(t, \cdot)$ under $\psi$.
Here $\mu=\frac{1}{2\pi} F_3$. 
One verifies that the map:
\begin{align*}
G:(h_-,h_+)\times \mu \mathbb{S} &\to A_{R_-, R_+} \subset \mathbb{C}\\
 (h, \theta)&\mapsto e^{\frac{h}{\mu}+i \frac{\theta}{\mu}}
 \end{align*} is a conformal diffeomorphism. 
 Here $A_{R_-, R_+}=\set{R_-<|z|<R_+}$ with $R_-=e^{\frac{h_-}{\mu}}$ and $R_+=e^{\frac{h_+}{\mu}}$.
As a consequence, we obtain a conformal diffeomorphism:
\begin{equation*}
\mathbf{F}=\psi \circ G^{-1}
\end{equation*} 
as in Lemma \ref{ShiffOssLem}.
We check that $\Sigma_t=\mathbf{F}( |z|=e^{ \frac{t}{\mu}} )$ and so
\begin{equation*}
\mathcal{H}^1(\Sigma_t)=L\left( \frac{t}{\mu}\right). 
\end{equation*}
The corollary then follows immediately from 
Lemma \ref{ShiffOssLem} and the fact that $\frac{1}{\mu^2}=\frac{(2\pi)^2}{F_3^2}$.
\end{proof}
\begin{rem}
We give an alternate approach to Corollary \ref{OSCor} in Appendix \ref{OtherApp}.  While this approach avoids the use of the Weierstrass representation and gives a sharper conclusion, it requires a certain geometric estimate that is still conjectural.
\end{rem}

\section{The Area Bound}
In order to prove Theorem $\ref{SecondThm}$ we use Corollary \ref{OSCor} to obtain a bound for the lengths of level sets:
\begin{prop} \label{LengthCompProp}
Let $P_-=\set{x_3=h_-}$ and $P_+=\set{x_3=h_+}$ be distinct parallel planes in $\Real^3$ with $h_-<h_+$ and let $\Omega$ be the open slab between them.   Fix $\Sigma \in \mathcal{A}(\Omega)$.  
Let $P_0=\set{x_3=h_0}\subset \bar{\Omega}$ denote the plane that satisfies:
\begin{equation*}
\mathcal{H}^1 (\Sigma \cap P_0)=\inf_{t\in (h_-,h_+)} \mathcal{H}^1(\Sigma_t).
\end{equation*}
Here $\Sigma_t=\Sigma\cap \set{x_3=t}$ and $\mathcal{H}^1 (\Sigma \cap P_+)$ is defined as $\liminf_{t\nearrow h_+ } \mathcal{H}^1 (\Sigma_t)$ and likewise for $\mathcal{H}^1 (\Sigma \cap P_-)$.
Let $C$ denote the vertical catenoid with $Flux(C)=(0,0,F_3(\Sigma))$, symmetric with respect to reflection through the plane $P_0$. 
If $C_t=C\cap \set{x_3=t}$ then for $t\in [h_-,h_+]$:
\begin{equation}
\mathcal{H}^1(\Sigma_t)\geq \mathcal{H}^1(C_t).
\end{equation}
 Equality can hold when $t\neq h_0$ if and only if  $\Sigma$ is a translate of $C\cap \Omega$.  
\end{prop}
\begin{proof}
Set $L_\Sigma(t)=\mathcal{H}^1(\Sigma_t)$ for $t\in (h_-,h_+)$.  By Lemma \ref{TransverseLem}, $L_\Sigma(t)$ depends smoothly on $t$ and by Lemma \ref{OSCor} one has:
\begin{equation*}
\frac{d^2}{dt^2} L_{\Sigma} (t)\geq \frac{(2\pi)^2}{F_3(\Sigma)^2} L_\Sigma (t),
\end{equation*}
with equality if and only if $\Sigma$ is piece of a catenoid.  
Notice that $L_\Sigma$ is a convex function on $(h_-,h_+)$.  By setting
$L_\Sigma(h_-)=\lim_{t\searrow h_- }L_\Sigma(t)$ and $L_\Sigma(h_+)=\lim_{t\nearrow h_+} L_\Sigma(t)$, we may think of $L_\Sigma$ as a function on $[h_-,h_+]$ but possibly taking the value $\infty$ at the end points.  The convexity ensures these limits exist.

For $C$ as in the statement of the theorem, $C\cap \Omega \in \mathcal{A}(\Omega)$.  Set $L_C(t)=\mathcal{H}^1(C_t)$.  As $Flux(C)\cdot \mathbf{e}_3=F_3(\Sigma)=F_3$ by assumption, Corollary \ref{OSCor} implies:
\begin{equation*}
\frac{d^2}{dt^2} L_{C} (t)= \frac{(2\pi)^2}{F_3^2} L_C (t).
\end{equation*}
Notice that $L_C$ is smooth on $[h_-,h_+]$ and the symmetry about $P_0=\set{x_3=h_0}$ implies $L'_C(h_0)=0$.

We claim that $L_\Sigma \geq L_C$ on $[h_-,h_+]$ with equality if and only if $\Sigma $ is a piece of a catenoid.  To see this we distinguish between when $h_0\in (h_-,h_+)$ and when $h_0$ is an endpoint. 
For any $t\in (h_-,h_+)$ 
\begin{equation*} 
F_3=\left|\int_{\Sigma_t} \mathbf{e}_3\cdot \mathbf{\nu} ds\right| \leq \int_{\Sigma_t} ds = L_\Sigma(t)
\end{equation*}
with equality if and only if $\Sigma_t$ is a geodesic in $\Sigma$. 
Similarly, if $t=h_-$ then 
\begin{equation*}
F_3=\lim_{t\searrow h_- } \left|\int_{\Sigma_t} \mathbf{e}_3\cdot \mathbf{\nu} ds\right|\leq \lim_{t\searrow h_- } \int_{\Sigma_t} ds = L_\Sigma(h_-)
\end{equation*}
 and the corresponding result holds when $t=h_+$. 
As $C_{h_0}$ is a geodesic in $C$,  $F_3=F_3(C)=L_C(h_0)$.  Hence, $L_C(h_0)\leq L_\Sigma(h_0)$. 

Now assume that $h_0\in (h_-,h_+)$.
For $t\in (h_-,h_+)$, the choice of $C$ and Corollary \ref{OSCor} ensure that  $L_C(h_0)\leq L_\Sigma(h_0)$, $L'_\Sigma(h_0)= 0=L'_C(h_0)$ 
and $\frac{d^2}{dt^2} \left( L_\Sigma(t)-L_C(t)\right) \geq \frac{(2\pi)^2}{F_3^2} \left( L_\Sigma(t)-L_C(t)\right)$. An ODE comparison 
then implies that for all  $t\in [h_-,h_+]$,  $L_C(t)\leq L_\Sigma(t)$ with equality holding for any $t\neq h_0$ if and only if $\Sigma$ is a piece of a catenoid and $\Sigma_{h_0}$ is a geodesic in $\Sigma$.   Thus, equality holds for any $t \neq h_0$, if and only if $\Sigma$ is equal (up to a translation) to $C\cap \Omega$.

When $h_0=h_-$ we argue as follows: For $\epsilon>0$ small, set $L_{C,\epsilon}(t)=L_C(t-\epsilon)$ and restrict attention to $(h_-+\epsilon, h_+)$.  Clearly, $L_{C,\epsilon}''=\frac{(2\pi)^2}{F_3^2} L_{C,\epsilon}$, $L'_{C,\epsilon}(h_-+\epsilon)=0$ and  $L_\Sigma(h_-+\epsilon)>L_\Sigma (h_-) \geq L_{C,\epsilon}(h_-+\epsilon)$.  Moreover, as $L_\Sigma$ is convex and has its minimum at $h_-$, $L_\Sigma'(h_-+\epsilon)>0$.  Hence by an ODE comparison,  $L_\Sigma(t)>L_{C,\epsilon}(t)$ for $t\in [h_-+\epsilon, h_+)$.  Letting $\epsilon\to 0$ implies $L_\Sigma(t)\geq L_{C}(t)$ for $t\in [h_-,h_+]$.  Equality can hold if and only if $\Sigma$ is equal (up to a translation) to $C\cap \Omega$.  An identical argument applies when $h_0=h_+$. 
\end{proof}
\begin{rem}
Proposition \ref{LengthCompProp} fails if $\Sigma$ were taken in the larger class of embedded elements of $\mathcal{M}(\Omega)$. 
Indeed, normalizing as in the proposition, it can be verified that (a suitably modified) version of Corollary \ref{OSCor} continues to hold for $\Sigma$ at $t$ so that $\set{x_3=t}$ meets $\Sigma$ transversally.  In particular, a modified version of Proposition \ref{LengthCompProp} holds between critical values of $x_3$.  However, it can be verified that while the length of level sets is continuous across critical values of $x_3$, the rate of change of the length of these level sets becomes infinite at a critical value.  In particular, there is never convexity across critical levels. 
\end{rem}

Let us now use Proposition \ref{LengthCompProp} to prove Theorem \ref{SecondThm}:
\begin{proof} 
We have verified:
\begin{equation*}
\mathcal{H}^1( C_t)\leq \mathcal{H}^1( \Sigma_t)
\end{equation*}
with equality for all $t$ if and only if $\Sigma$ is equal (up to a translation) to $C\cap \Omega$. Fix $h_*\in(h_-,h_+)$ and define the function $A_{\Sigma,h_*}(t)$ on $[h_*, h_+) $ by 
 \begin{equation*}
 A_{\Sigma,h_*} (t)=\mathcal{H}^2(\Sigma \cap \set{h_*\leq x_3 < t})
 \end{equation*}
 and the function $A_{C,h_*}$ similarly. 

The co-area formula implies that
\begin{equation*}
\frac{d}{dt}A_{\Sigma,h_*}(t)= \int_{\set{x_3=t}\cap \Sigma} \frac{1}{|\nabla_\Sigma x_3|} .
\end{equation*}
Applying the Cauchy-Schwarz inequality yields:
\begin{equation*}
\frac{d}{dt}A_{\Sigma,h_*}(t)\geq \frac{\mathcal{H}^1 (\Sigma_t)^2}{\int_{\Sigma_t} |\nabla_\Sigma x_3| }=\frac{\mathcal{H}^1(\Sigma_t)^2}{F_3}.
\end{equation*}
Notice that one has equality if and only if $ \frac{1}{|\nabla_\Sigma x_3|}$ and $|\nabla_\Sigma x_3|$ are linearly dependent, in other words are both constant.  This is readily checked to be the case on $C$ and so using the above estimate for length:
\begin{equation*}
\frac{d}{dt}A_{\Sigma,h_*}(t)\geq \frac{\mathcal{H}^1( C_t)^2}{F_3 }=\frac{d}{dt}A_{C,h_*}(t).
\end{equation*}
Integrating implies $A_{\Sigma,h_*} (t)\geq A_{C, h_*}(t)$ for  $t\in [h_*,h_+)$ with equality if and only if $\Sigma$ is a piece of of $C$.  Letting $h_*\to h_-$ proves the theorem.
\end{proof}

\section{Sharp non-existence result}
\label{NoSpanSec}
As discussed in the introduction, given an open slab $\Omega$ with $\partial \Omega=P_+\cup P_-$ and connected, simple closed curves $\sigma_\pm \subset P_\pm$ there need not be a surface $\Sigma\in \mathcal{M}(\Omega)$ with $\partial \Sigma=\sigma_+\cup \sigma_-$.  For instance, if the curves $\sigma_\pm$ are too short relative to the height of the slab then there cannot be a connected minimal surface $\Sigma\in \mathcal{M}(\Omega)$ spanning $\sigma_\pm$.  Indeed, the monotonicity formula gives a lower bound on the area of such a $\Sigma$ in terms of the distance between the planes, while the isoperimetric inequality gives an upper bound in terms of the lengths of the curves  (for surfaces in $\mathcal{M}(\Omega)$ with two boundary components the isoperimetric inequality with sharp constant is known to hold -- see \cite{LiSchoenYau}). Alternatively, if the $\sigma_{\pm}$ are well separated, barrier arguments can be used to rule out the existence of such $\Sigma$.   Using Proposition \ref{LengthCompProp}, we are able to  give a sharp condition (see also Theorem 6 of \cite{OS} for a related result):
\begin{thm}\label{NoSpanThm}
Fix $\Omega$ an open slab with $\partial \Omega=P_-\cup P_+$ the union of two parallel planes.  Let $\sigma_\pm \subset P_\pm$ be a pair of connected simple closed curves. Let $C_{MS}$ be the unique (up to translations parallel to $P_\pm$) minimal surface in $\mathcal{A}(\Omega)$ obtained via rigid motions and homotheties from $Cat_{MS}$.  If we define $L_{crit}(\Omega):=\mathcal{H}^1(\partial C_{MS})$ and
\begin{equation*}
 \mathcal{H}^1(\sigma_+\cup \sigma_-)<L_{crit}(\Omega)
\end{equation*}
 then there is no surface $\Sigma\in \mathcal{M}(\Omega)$ with $\partial \Sigma=\sigma_+\cup \sigma_-$.  Moreover, if $\Sigma \in \mathcal{M}(\Omega)$ is a smooth minimal surface with $\partial \Sigma =\sigma_-\cup \sigma_+$ and
\begin{equation*}
 \mathcal{H}^1(\partial \Sigma)=L_{crit}(\Omega)
\end{equation*}
 then $\Sigma$ is a translate of $C_{MS}$.
\end{thm}

In order to prove this theorem, we first prove a more general result.  Namely, we will show the existence of a 
$\Sigma \in \mathcal{M}(\Omega)$ with $\partial \Sigma = \sigma_-\cup \sigma_+$ is precluded if one boundary curve is 
too short as determined by an explicit function of the length of the other boundary curve.  Roughly speaking, the existence of
such a $\Sigma$ relies on the existence of a vertical catenoid $C$ so that $\mathcal{H}^1(C\cap P_\pm)= \mathcal{H}^1(\sigma_+ \cup \sigma_-)$. We point out that Theorem 4 of \cite{OS} gives the same result when one considers only $\Sigma\in \mathcal{A}(\Omega)$. As in the case 
for area bounds, a marginally stable piece of a catenoid will serve as the model.  However, here the marginally stable pieces are 
generally not obtained from rigid motions and homotheties of $Cat_{MS}$.

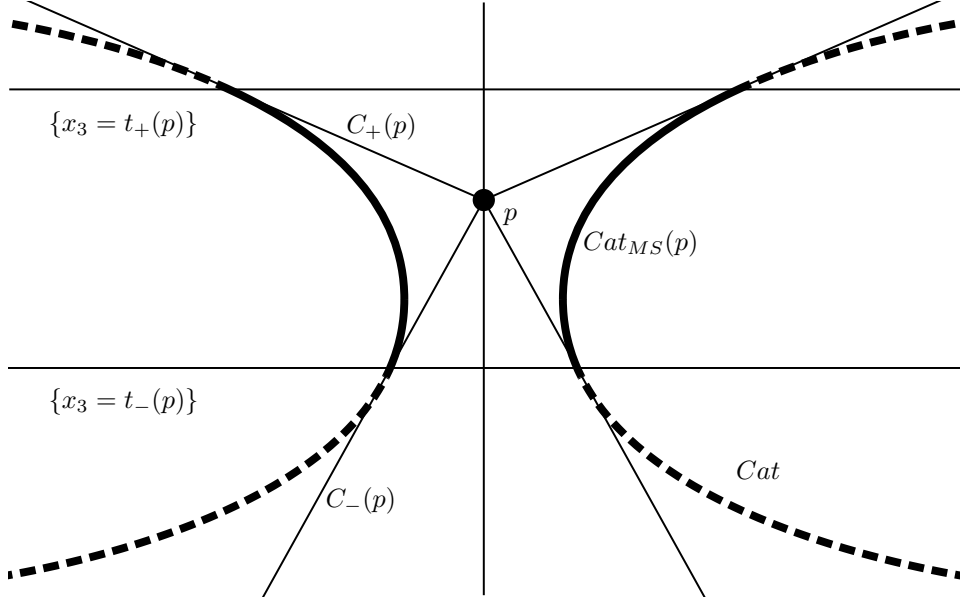
\begin{figure}
 \centering
\def\svgwidth{\columnwidth}
\input{MSCat2.tex}
\caption{The subset $Cat_{MS}(p)$ of $Cat$ is indicated as well as the cones $C_+(p)$ and $C_-(p)$. }
\label{MSCatFig}
\end{figure}

We begin by describing the general class of marginally stable pieces of $Cat$ we will need.
First note that the 
rotational symmetry and convexity of the function $\cosh t$ imply that for each point $p=z\mathbf{e}_3$ on the $x_3$-axis, there are unique cones over $p$ that intersect $Cat$ tangentially.  Precisely, there exist values $t_+=t_+(p)>0$ and $t_-=t_-(p)<0$
 with the following property: the cones $C_+(p)$ (resp. $C_-(p)$) over $p$ of $Cat_{t_+(p)}=Cat\cap\set{x_3=t_+(p)}$ (resp. 
$Cat_{t_-(p)}=Cat\cap\set{x_3=t_-(p)}$) meet $Cat$ only at $Cat_{t_+(p)}$ (resp. $Cat_{t_-(p)}$).  We observe also that $t_{+}$ is an increasing and continuous 
function of $z$ with range $(0,\infty)$; similarly, $t_-$ is increasing and continuous with range $(-\infty,0)$.  Notice that $Cat$  must be tangential to $C_\pm(p)$ at $Cat_{t_\pm(p)}$. We refer the reader to Figure \ref{MSCatFig}.

 Let $Cat_{MS}(p)$ be the bounded component of $Cat\backslash \left(C_+(p)\cup C_-(p)\right)$.  
One verifies that $Cat_{MS}=Cat_{MS}(0)$ and that as $p\to (0,0,\infty)$, $Cat_{MS}(p)$ 
converges to $Cat\cap \set{x_3>0}$.  We claim that $Cat_{MS}(p)$ is marginally stable for each $p$.  This 
follows from  the observation that for $\lambda>0$ the surfaces $Cat_{MS}^\lambda(p)=p+\lambda\left(Cat_{MS}(p)-p\right)$  
give a foliation of the component $\bar{C}(p)$ of $\Real^3\backslash \left(C_{+}(p)\cup C_-(p)\right)$ containing $Cat_{MS}(p)$.
 Moreover,  as $Cat$ meets $C(p)=\partial \bar{C}(p)$ tangentially,  the motion of $\partial Cat_{MS}^\lambda(p)$ at $\lambda=1$ is 
tangential to $Cat_{MS}(p)$.  As a consequence, the normal variation at $\lambda=1$ gives a positive 
Jacobi field on $Cat_{MS}(p)$ vanishing on $\partial Cat_{MS}(p)$.  One also verifies that if $C_0=Cat\cap \set{h_-<x_3<h_+}$ is 
marginally stable then $C_0=Cat_{MS}(p)$ for exactly one value $p$.  

 Using these marginally stable pieces, we now determine the explicit function of interest:
 \begin{lem} \label{MSExistLem}
Fix $\Omega_0=\set{-1<x_3<1}$ and $\partial \Omega_0=P_-\cup P_+$.  There exist two well-defined functions $C_{MS}:\Real^+ \to \mathcal{A}(\Omega_0)$ and $F_{\Omega_0}:\Real^+\to \Real^+$ determined in the following way: For each $L_->0$ let $C_{MS}(L_-)$ denote the unique (up to translations parallel to $P_\pm$) marginally stable piece of a vertical catenoid with 
\begin{equation*}
 \partial C_{MS}(L_-)=\sigma_-\cup \sigma_+\subset P_-\cup P_+
\end{equation*}
and $\mathcal{H}^1(\sigma_-) = L_-$. Given $C_{MS}(L_-)$, define
\begin{equation*}
 F_{\Omega_0}(L_-)=\mathcal{H}^1(\sigma_+).
\end{equation*}
Furthermore, $F_{\Omega_0}$ has the following properties:
\begin{enumerate}
\item \label{PropOne}If $L_+\geq F_{\Omega_0}(L_-)$ then there is a vertical catenoid $C=C(L_-,L_+)$ so that writing $\partial(C \cap \Omega_0)=\gamma_-\cup \gamma_+$ gives $\mathcal{H}^1(\gamma_\pm)=L_\pm$, while if $L_+<F_{\Omega_0}(L_-)$ no such vertical catenoid exists.
\item \label{PropTwo}If $C$ is a vertical catenoid with $\partial(C\cap \Omega_0)=\gamma_-\cup \gamma_+$ and $\mathcal{H}^1(\gamma_-)<L_-$ then $\mathcal{H}^1(\gamma_+)> F_{\Omega_0}(L_-)$.
\item \label{PropThree} If $C$ is a vertical catenoid with $\partial(C\cap \Omega_0)=\gamma_-\cup \gamma_+$, $\mathcal{H}^1(\gamma_-)=L_-$, and $\mathcal{H}^1(\gamma_+)=F_{\Omega_0}(L_-)$ then $C\cap \Omega_0$ is a translate of $C_{MS}(L_-)$.
\end{enumerate}
\end{lem}
\begin{rem}
We note that for other slabs $\Omega$ it is straightforward to determine $F_{\Omega}$ in terms of $F_{\Omega_0}$. Indeed, rigid motions leave the function invariant and $F_{\lambda \Omega}(L)=\lambda F_{\Omega}(\lambda^{-1} L)$.
\end{rem}
\begin{proof}
We claim that given any $L_->0$ there is a marginally stable piece of a vertical catenoid $C_{MS}(L_-)$ with $\partial C_{MS}(L_-)=\sigma_- \cup \sigma_+\subset P_-\cup P_+$ and $\mathcal{H}^1(\sigma_-)=L_-$.  Rather than prove this by direct computation, we use global arguments.  Set
\begin{equation*}
C=\frac{L_-}{2\pi} Cat -\eE_1.
\end{equation*} 
Let $C_0=C\cap \Omega_0$ and denote $\partial C_0=\sigma_-^0\cup \sigma_+^0$ so that  $\mathcal{H}^{1} (\sigma_-^0)=L_-$.  By domain monotonicity for eigenvalues, $C_0$ is strictly stable because $C\cap \set{x_3\geq -1}$ is stable.  
 Consider now the following smooth family of coaxial circles in $P_\pm$: for $t>0$ set  $\sigma_-^t=\sigma_-^0$ and  $\sigma_+^t=t(\sigma_+^0-\eE_1)+\eE_1$.  By \cite{MeeksWhite} and a barrier argument, there is a $1>T_{crit}>0$ so that for each $t\in (T_{crit},1]$ there are $C_t$, strictly stable minimal annuli smoothly depending on $t$, with $\partial C_t=\sigma_-^t\cup \sigma_+^t$ and for $t=T_{crit}$ there is a marginally stable annulus, $C_{T_{crit}}$, with $\partial C_{T_{crit}}=\sigma_-^{T_{crit}}\cup \sigma_+^{T_{crit}}$. As the boundaries consist of coaxial circles, the proof of \cite{SCH} implies that each $C_{t}$ is a piece of a catenoid.    
 In fact, for all $t\in (T_{crit},\infty)$ there is a strictly stable minimal annulus $C_t$ with $\partial C_t=\sigma_-^t\cup \sigma_+^t$.  To verify the claim, it suffices to consider $t\in (1,\infty)$ and in this range the $C_t$ are obtained from appropriate rescalings and translations of subsets of $C$.

We claim that $C_{T_{crit}}$ is the desired $C_{MS}(L_-)$ and gives $F_{\Omega_0}$ as outlined.  
We first note that by the uniqueness of the sets $Cat_{MS}(p)$ there are $p,\lambda$ and $h$ so that $C_{T_{crit}}=\lambda (Cat_{MS}(p)-p)+p+h\mathbf{e}_3$. That is, $C_{T_{crit}}$ is, up to a vertical translation, a rescaling of one of the marginally stable pieces described previously.

For $C'$ a vertical catenoid, let $C_0'=C'\cap \Omega_0$. Suppose $C'$ is such that $\partial C_0'=\sigma_-'\cup \sigma_+'\subset P_-\cup 
P_+$ and $\mathcal{H}^1(\sigma_-')=L_-$. For the same $p$ as before, we may write $C_0'=(\lambda' (Cat-p)+p +(h+\Delta h)\mathbf{e}_3)\cap \Omega_0$. We claim $\Delta h \leq 0$.  This follows 
by noting that as $Cat_{MS}(p)$ meets $C(p)$ tangentially $\left(\lambda' (Cat-p)+p+h\mathbf{e}_3\right) \cap P_-$ has length greater than $L_-$. An 
upward translation only increases the length further and so one must translate downward, i.e. take $\Delta h\leq 0$. By a similar reasoning, 
one concludes that $\mathcal{H}^1(\sigma'_+)\geq \mathcal{H}^1(\sigma_+^{T_{crit}})$.  As $\sigma_\pm'$ are coaxial circles, 
we have 
$\sigma_-'=\sigma_-^{T_{crit}}$ and $\sigma_+'$ surrounding $\sigma_+^{T_{crit}}$. By \cite{MeeksWhite}, as $C_{T_{crit}}$ is 
marginally stable, if $C_0'$ differs from $C_{T_{crit}}$ then $\partial C_0'$ differs from $\partial C_{T_{crit}}$; hence in this case 
$\sigma_+'$ strictly surrounds $\sigma_+^{T_{crit}}$.  By the preceeding paragraph we then see that $\partial C_0'$ bounds a strictly stable
 annulus and so, by \cite{MeeksWhite}, $C_0'$ cannot be marginally stable.  This proves the claimed uniqueness. 

 Clearly, \eqref{PropOne} is an immediate consequence of the preceeding argument.  Furthermore, if \eqref{PropTwo} failed to hold for a vertical catenoid $C'$ then it could be used as a barrier allowing one to construct a piece of a vertical catenoid  violating \eqref{PropOne} (see the proof of Proposition \ref{SharpCurvProp} for a detailed argument).  Finally, by \cite{MeeksWhite}, as $C_{T_{crit}}$ is marginally stable there is no other minimal surface $C_{T_{crit}}'$ with $\partial C_{T_{crit}}=\partial C_{T_{crit}}'$ which verifies \eqref{PropThree}.
\end{proof}

As a consequence we may prove the following general proposition giving sharp conditions for the non-existence of minimal surfaces spanning a given pair of curves:
\begin{prop} \label{SharpCurvProp}
Fix $\Omega$ an open slab with $\partial \Omega=P_-\cup P_+$ the union of two parallel planes and let $C_{MS}:\Real^+ \to \mathcal{A}(\Omega), \;F_{\Omega}:\Real^+\to \Real^+$ be the functions given by Lemma \ref{MSExistLem}.  Let $\sigma_\pm \subset P_\pm$ be a pair of connected, simple closed curves. If 
\begin{equation*}
\mathcal{H}^1(\sigma_+)<F_\Omega(\mathcal{H}^1(\sigma_-))
\end{equation*}
then there is no surface $\Sigma\in \mathcal{M}(\Omega)$ with $\partial \Sigma=\sigma_+\cup \sigma_-$.  Moreover, if $\Sigma\in \mathcal{M}(\Omega)$ is a smooth minimal surface with $\partial \Sigma =\sigma_-\cup \sigma_+$ such that
\begin{equation*}
\mathcal{H}^1(\sigma_+)=F_\Omega(\mathcal{H}^1(\sigma_-))
\end{equation*}
  then  $\Sigma$ is a translate of $C_{MS}(\mathcal{H}^1(\sigma_-))$.
\end{prop}
\begin{rem}  Let $\mathcal{M}_2(\Omega)\subset \mathcal{M}(\Omega)$ be the set of all $\Sigma \in \mathcal{M}(\Omega)$ with $b\Sigma = \partial \Sigma$ consisting of exactly two connected boundary components.
 If one considers $\Psi_\Omega: \mathcal{M}_2(\Omega)\to \Real^2$ the map defined by 
\begin{equation*}
\Psi(\Sigma)=\left( \mathcal{H}^1(\partial \Sigma\cap P_-), \mathcal{H}^1(\partial \Sigma\cap P_+)\right)
\end{equation*} 
then the proposition says that the image of $\Psi$ is an unbounded region in the first quadrant of the plane whose boundary consists of the images of marginally stable pieces of catenoids and is explicitly given as the graph of the function $F_{\Omega}$.  
\end{rem}

\begin{proof}
Up to a rescaling and rigid motion we may take $\Omega=\set{-1<x_3<1}$.
 Suppose that $\Sigma\in \mathcal{M}(\Omega)$ has the structure of a smooth manifold with boundary and that $\partial \Sigma$ is embedded and consists of two connected components $\sigma_\pm \subset P_\pm$. By assumption, the $\sigma_\pm$ are connected, simple closed curves in $P_+$ and $P_-$.  It will suffice to show that $\mathcal{H}^1 (\sigma_+)\geq F_{\Omega}(\mathcal{H}^1(\sigma_-))$.   Note that $\Sigma$ is allowed to be immersed and have arbitrary genus, however it may still be used as a barrier to construct an embedded annulus with the same boundary.  Indeed, while $\Omega\backslash \Sigma$ may have more than $2$ components, only one of these, $\Omega'$, is unbounded.  Clearly, $\sigma_+$ and $\sigma_-$ are homotopic in $\bar{\Omega}'$ but are not null homotopic in $\bar{\Omega'}$.  In particular, there is an annulus $A$ in $\Omega'$ with $\partial A=\sigma_+\cup \sigma_-$ but no disk $D$ in $\Omega'$ with $\partial D=\sigma_+$ or $\partial D=\sigma_-$   Finally, we point out that $\bar{\Omega}'$ is mean convex in the sense of Meeks and Yau \cite{MeeksYau2}.  As a consequence, by \cite{MeeksYau2} there is an embedded minimal annulus $\Gamma\subset \Omega'$ with $\partial \Gamma =\sigma_+\cup \sigma_-$.

By Proposition \ref{LengthCompProp}, there is a vertical catenoid $C$ so that if we write $\partial(C\cap \Omega)=\gamma_+\cup \gamma_-$ where $\gamma_\pm\subset P_\pm$ then $\mathcal{H}^1(\sigma_\pm)\geq \mathcal{H}^1(\gamma_\pm)$. Moreover the inequality is strict unless $C\cap \Omega$ and $\Gamma$ agree up to a translation parallel to $P_\pm$.  As $\mathcal{H}^1(\gamma_-)\leq \mathcal{H}^1(\sigma_-)$, by \eqref{PropTwo} of Lemma \ref{MSExistLem} 
\begin{equation*}
\mathcal{H}^1(\sigma_+)\geq \mathcal{H}^1(\gamma_+)\geq F_{\Omega}
( \mathcal{H}^1(\sigma_-)).
\end{equation*}
Finally, by Lemma \ref{MSExistLem} and Proposition \ref{LengthCompProp} equality is only achieved if $\Gamma$ is a horizontal translate of $C_{MS}(\mathcal{H}^1(\sigma_-))$.  In this case,  $\sigma_+\cup \sigma_-=\partial \Gamma$ consist of coaxial circles in parallel planes.  Hence, the proof of \cite{SCH} implies the $\sigma_\pm$ can bound only pieces of a catenoid; that is, $\Sigma$ is a horizontal translate of $C_{MS}(\mathcal{H}^1(\sigma_-))$.
\end{proof}

We now prove Theorem \ref{NoSpanThm}:
\begin{proof}Up to a rescaling and rigid motion we may take $\Omega=\set{-1<x_3<1}$.
By Proposition \ref{SharpCurvProp} we need only verify the theorem for vertical catenoids. 
  The space of vertical catenoids is parameterized by $\lambda>0$ and $t$ where
\begin{equation*}
Cat_{\lambda,t}=\lambda Cat+t\mathbf{e}_3.
\end{equation*}
Just as in the proof of Theorem \ref{MainThm} where we saw $C_{MS}$ minimized area among all vertical catenoid pieces in $\Omega$, we now show $C_{MS}$ minimizes boundary length in this same class. Setting $L(\lambda,t)=\mathcal{H}^1(\partial(\Omega\cap Cat_{\lambda,t}))$,
one computes 
\begin{equation*}
L(\lambda,t)=2\pi \lambda \cosh \frac{1-t}{\lambda}+2\pi \lambda \cosh \frac{-1-t}{\lambda}.
\end{equation*}
As
$L(\lambda,t)\to \infty $ when $|(\ln \lambda,t)|\to \infty$, it suffices to find critical points of $L$.  First observe that $\frac{\partial}{\partial t} L=-2\pi \sinh \frac{1-t}{\lambda}-2\pi  \sinh \frac{-1-t}{\lambda}$, which is zero only when $t=0$.  Thus one must only minimize $L(\lambda,0)=4\pi \lambda \cosh \frac{1}{\lambda}$.  One verifies that the critical points of $L(\lambda,t)$ are of the form $(\lambda_0,0)$ where $\lambda_0$ satisfes $\lambda_0=\tanh \frac{1}{\lambda_0}$. Hence, as in Theorem \ref{MainThm}, $\lambda_0$ is unique and $Cat_{\lambda_0,0}=C_{MS}$.
\end{proof}

\appendix
\section{Two Lemmas of Osserman and Schiffer} \label{OSApp}
For the sake of completeness, we present here a proof of Lemma \ref{ShiffOssLem}. The argument is that given by Osserman and Schiffer in \cite{OS} though we have updated the notation where appropriate and omitted some details.   The argument makes crucial use of the Weierstrass representation and so we first discuss this fundamental connection between minimal surfaces in $\Real^3$ and complex analysis.

Consider as in Lemma \ref{ShiffOssLem} the following conformal, harmonic immersion:
\begin{equation*}
\mathbf{F}: A_{r, R}\to \Real^3.
\end{equation*}
Here $A_{r,R}=\set{r<|z|<R}\subset \mathbb{C}$.
In particular, the image $\Sigma$ of $\mathbf{F}$ is minimal.
Denote by $dh=h(z) dz$ the holomorphic one form on $A_{r,R}$ whose real part is $\mathbf{F}^{*} dx_3$ and by $g$ the function on $A_{r,R}$ given by the stereographic projection of the normal of $\Sigma$.  Here $z$ is the coordinate on $A_{r,R}$ induced from $\mathbb{C}$.  It is a standard exercise to see that the minimality of the image $\Sigma$ implies that $g$ is meromorphic. The Weierstrass representation allows one to recover $\mathbf{F}$ (up to a translation) from the data $dh$ and $g$.
Indeed, one has:
\begin{equation}\label{WeierstrassRep}
\mathbf{F}:=\re \int\left(\frac{1}{2}
(g^{-1}-g),\frac{i}{2}(g^{-1}+g), 1\right) dh.
\end{equation}
We point out that $g$ and $dh$ are not arbitrary.  
Indeed, it is straightforward to compute:
\begin{equation} \label{WMetricEqn}
\mathbf{F}^*g_{euc}=\frac{1}{4} (|g|+|g|^{-1})^2 dh \otimes d \bar{h}=\frac{|h|^2}{4} (|g|+|g|^{-1})^2 dz \otimes d\bar{z}.
\end{equation}
Hence, as $\mathbf{F}$ is an immersion, $hg$ and $hg^{-1}$ are both holomorphic and cannot vanish simultaneously. 
Additionally, as closed curves in $A_{r,R}$ should map to closed curves in $\Sigma$, the Weierstrass data must satisfy the following \emph{period conditions}:
\begin{equation*}
\int_{\gamma} g dh =\overline{\int_{\gamma} g^{-1} dh}, \; \; \; \re \int_{\gamma} dh =0
\end{equation*}
for any closed curve $\gamma$ in $A_{r,R}$. 
Finally, the Weierstrass data can be used to compute $Flux(\mathbf{F}(\gamma))$ for a curve $\gamma\subset A_{r,R}$:
\begin{equation}\label{FluxDefEq}
Flux(\mathbf{F}(\gamma))=\im \int_{\gamma} \left( \frac{1}{2}
(g^{-1}-g),\frac{i}{2}(g^{-1}+g), 1\right) dh.
\end{equation}

Lemma \ref{ShiffOssLem} is a simple consequence of the following (also reproduced from \cite{OS}):
\begin{lem}\label{CPXLem}
Let $F$ be a holomorphic function on the annulus $A_{r,R}$.  If $F$ has no zeros on $\set{|z|=\rho}\subset A_{r,R}$ and satisfies $\int_{|z|=\rho} F(z) \frac{dz}{z}=0$ then
\begin{equation*}
\int_{|z|=\rho} \rho^2 \Delta |F|  \left| \frac{dz}{z} \right| \geq \int_{|z|=\rho} |F|   \left| \frac{dz}{z} \right|
\end{equation*}
with equality if and only if $F=a z$ or $F=a z^{-1}$.
\end{lem}
\begin{proof}
Let $G$ be an arbitrary holomorphic function on the annulus $A_{r',R'}$ and set $a_0=\int_{|z|=\rho} G  \left| \frac{dz}{z} \right|$.
As $G$ is holomorphic:
\begin{equation*}
\Delta |G|^2= 4|G'|^2.
\end{equation*}
Thus, the Cauchy-Riemann equations and the Wirtinger inequality imply:
\begin{equation} \label{GIneq}
\int_{|z|=\rho} \rho^2 \Delta |G|^2  \left| \frac{dz}{z} \right| \geq 4 \int_{|z|=\rho} |G|^2  \left| \frac{dz}{z} \right| -8 \pi |a_0|^2.
\end{equation}

As $F$ is non-vanishing on $|z|=\rho$ there are $r<r'<R' <R$ so that $r'<\rho<R'$ and $F$ is non-vanishing on $A_{r',R'}$.  In particular, the winding number of the map $F:A_{r',R'}\to \mathbb{C}\backslash \set{0}$ is a well-defined integer $k$ with
\begin{equation*}
\int_{|z|=\rho} \frac{dF}{F}= 2\pi i k.
\end{equation*}
If $k$ is even then there is holomorphic function $G$ on $A_{r',R'}$ so that $F=G^2$, while if $k$ is odd there is holomorphic function $G$ on $A_{r',R'}$ so that $F=z G^2$.  In both cases expand $G$ in a Laurent series as
\begin{equation}\label{LExpansion}
G(z)=\sum_{n=-\infty}^\infty a_n z^n.
\end{equation}
We treat the two cases separately:

{\bf Case 1}:  When $k$ is even,  the constant term in the Laurent expansion of $F$ is given (in terms of \eqref{LExpansion}) by:
\begin{equation*}
a_0^2+2\sum_{n=1}^\infty a_{n} a_{-n}=0;
\end{equation*}
 here the  condition that $F\frac{dz}{z}$ integrate to $0$ on $|z|=\rho$ is used. 
Hence,
\begin{align*}
|a_0|^2 &=2\left|\sum_{n=1}^\infty a_{n}\rho^{n} a_{-n}\rho^{-n}\right|\\
 &\le	 \sum_{n=1}^\infty \left( a_{n}^2 \rho^{2 n} + a_{-n}^2 \rho^{-2 n}\right)\\
 &=\frac{1}{2\pi} \int_{|z|=\rho} |G|^2 \left| \frac{dz}{z} \right| -|a_0|^2.
 \end{align*}
 That is,
 \begin{equation*}
 4\pi |a_0|^2\leq \int_{|z|=\rho} |G|^2  \left| \frac{dz}{z} \right| .
 \end{equation*}
 Combining this with \eqref{GIneq} proves the Lemma in this case.

 {\bf Case 2}: As $k$ is odd we have $F=z G^2$.  We introduce an auxilliary holomorphic function $H=z  G$ so that $z F=H^2$.  One computes:
 \begin{align*}
 |H'|^2= \frac{|z F' +F|^2}{4|zF|}, \\
  |G'|^2=\frac{|z F'-F|^2}{4|z^3 F|}.
 \end{align*}
 The constant term in the Laurent expansion of $H$ is the term $a_{-1}$ from \eqref{LExpansion}.
 Hence, the Wirtinger inequality gives:
 \begin{align*}
 \int_{|z|=\rho} |H|^2 \left| \frac{dz}{z} \right| &\leq   \rho^2 \int_{|z|=\rho}  |H'|^2  \left| \frac{dz}{z} \right| +2 \pi |a_{-1}|^2\\
 &= \frac{\rho^2}{4}\int_{|z|=\rho}\frac{|z F' +F|^2}{|zF|} \left| \frac{dz}{z} \right|+2\pi |a_{-1}|^2. \\
  \end{align*}
Similarly,
\begin{align*}
\int_{|z|=\rho} |G|^2 \left| \frac{dz}{z} \right|&\leq \rho^2 \int_{|z|=\rho} |G'|^2 \left| \frac{dz}{z} \right|+2\pi |a_0|^2\\
&\leq  \frac{\rho^2}{4}\int_{|z|=\rho}  \frac{|z F' -F|^2}{|z^3 F|} \left| \frac{dz}{z} \right|+2\pi |a_0|^2. 
\end{align*}
As $|zF'+F|^2+|zF'-F|^2 =2|z F'|^2 +2|F|^2$, combining the two inequalities and using $G^2=\frac{F}{z}, H^2=z F$ yields:
\begin{equation*}
3 \int_{|z|=\rho} |F| \left| \frac{dz}{z} \right|\leq \int_{|z|=\rho} \rho^2 \frac{|F'|^2}{|F|} \left| \frac{dz}{z} \right|+4\pi \rho \left(\frac{|a_{-1}|^2}{\rho^2}+|a_0|^2\right).
\end{equation*}
However,
\begin{equation*}
\int_{|z|=\rho} |F| \left| \frac{dz}{z} \right|=\rho \int_{|z|=\rho} |G|^2 \left| \frac{dz}{z} \right|\geq 2\pi \rho \sum_{n=-\infty}^\infty |a_n|^2 \rho^{2n} \geq 2\pi \rho \left( \frac{|a_{-1}|^2}{\rho^2} +|a_0|^2 \right)
\end{equation*}
and thus
\begin{align*}
\int_{|z|=\rho} |F| \left| \frac{dz}{z} \right|&\leq \int_{|z|=\rho} \rho^2 \frac{|F'|^2}{|F|} \left| \frac{dz}{z} \right|= \int_{|z|=\rho} \rho^2 \Delta |F| \left| \frac{dz}{z} \right|
\end{align*}
where the last equality follows from the fact that $F$ is holomorphic.

To see that the equality holds if and only if $F=az $ or $F=\frac{a}{z}$ we refer the reader to \cite{OS}.
\end{proof}

We now prove Lemma \ref{ShiffOssLem}:
\begin{proof}
Let $dh$ and $g$ be the Weierstrass data associated to $\mathbf{F}$.  Let $\Sigma$ denote the minimal annulus that is the image of $\mathbf{F}$.  For $z$ the usual coordinate on $A_{r,R}$,  $g=g(z)$ is a meromorphic function and $dh=h(z) dz$ for some holomorphic function $h(z)$.
Up to taking an ambient rotation of $\Real^3$ we may assume that $Flux(\Sigma)=(0,0,\lambda)$, that is $\Sigma$ has vertical flux.  
Using \eqref{FluxDefEq}, it is straightforward to verify this flux condition restricts the Weierstrass data as follows:
\begin{equation*}
\int_{\gamma} g dh =-\overline{\int_{\gamma} g^{-1} dh}.
\end{equation*}
On the other hand the period conditions imply:
\begin{equation*}
\int_{\gamma} g dh =\overline{\int_{\gamma} g^{-1} dh}
\end{equation*}
and so one concludes:
\begin{equation*}
\int_{\gamma} g dh=\int_{\gamma} zgh \frac{dz}{z}=0, \int_{\gamma} g^{-1} dh =\int_{\gamma} z hg^{-1} \frac{dz}{z}=0.
\end{equation*}
If $\sigma_\rho$ is the curve $\mathbf{F}(\set{|z|=\rho})$ we compute using \eqref{WMetricEqn} and the area formula:
\begin{align*}
\mathcal{H}^1(\sigma_\rho)&=\frac{1}{2}\int_{|z|=\rho}  \left(|z g h|+|z h g^{-1}|\right)\frac{|dz|}{|z|}.
\end{align*}
Recall, $zgh$ and $zhg^{-1}$ are holomorphic and non-vanishing, and so have zeros only at a discrete set of points.  We assume in what follows that neither have a zero on $|z|=\rho$.

If we denote by $\bar{L}(\rho)=\mathcal{H}^1(\sigma_\rho)$ then $L(t)=\bar{L}(e^t)$.  Thus,
\begin{align*}
L''(t)&=\left.\rho \frac{d}{d\rho}\left( \rho \frac{d}{d\rho} \bar{L}(\rho)\right) \right|_{\rho=e^t}\\
	&=\left.\frac{1}{2}\int_{|z|=\rho}  \rho \frac{d}{d\rho}\left( \rho \frac{d}{d\rho} \left(|z g h|+|z h g^{-1}|\right) \right) \left| \frac{dz}{z} \right| \right|_{\rho=e^t}\\
	&=\left. \frac{1}{2}\int_{|z|=\rho}  \rho^2 \Delta \left(|z g h|+|z h g^{-1}|\right)  \left| \frac{dz}{z} \right| \right|_{\rho=e^t}
\end{align*}
where the last equality follows from integration by parts and the fact that on $A_{r,R}$:
\begin{equation*}
\Delta =\frac{1}{\rho} \frac{d}{d\rho}\left( \rho \frac{d}{d\rho}  \right) +\frac{1}{\rho^2} \frac{d^2}{d\theta^2}.
\end{equation*}
   
Notice that we have already ensured that
\begin{equation*}
\int_{|z|=\rho}  z hg \frac{dz}{z} =\int_{|z|=\rho} zh g^{-1} \frac{dz}{z} =0
\end{equation*}
and so Lemma \ref{CPXLem} proves the theorem provided $zhg$ and $zhg^{-1}$ do not have a zero on $|z|=e^t$.  Hence, away from a finite number of $t$ the inequality holds.   However, as $zhg$ and $z hg^{-1}$ cannot simultaneously vanish $L''(t)$ is continuous in $t$.  Indeed, this is clear once one notes:
\begin{equation*}
 |zhg|+|zhg^{-1}|=|z|\sqrt{|hg|^2 +| hg^{-1}|^2+2|h|^2}
\end{equation*}
and that $|z|>0$.  As $L(t)$ is then also continuous, the result holds for all $t$.    For the case of equality we refer to \cite{OS}.
\end{proof}
\section{A conjectural approach}\label{OtherApp}
One downside to the use of Lemma \ref{ShiffOssLem} is that it depends in an essential manner on the Weierstrass representation.  This has the disadvantage of obscuring some of the geometric meaning as well as restricting applications to minimal surfaces in $\Real^3$.
For both these reasons it is fruitful to find a proof that avoids the use of the Weierstrass representation. 
In this section we give such an approach, albeit with one important caveat.  Namely, we require a certain sharp eigenvalue estimate to hold that is, to our knowledge, still conjectural.   We feel justified in presenting this approach both for the reasons already mentioned and because the conjecture is geometrically natural and seems to have broader applications in spectral theory.

\begin{conj} \label{OvalConj}
Let $\sigma$ be a smooth closed curve in $\Real^3$ parameterized by arclength $s$.  Denote by $\kappa$ the geodesic curvature of $\sigma$. 
Then for any smooth function $f$ on $\sigma$
\begin{equation} \label{BestConstIneq}
\mathcal{H}^1(\sigma)^2 \int_\sigma \left(\left|\frac{d f}{ds}\right|^2 +\kappa^2 f^2 \right)ds \geq (2\pi)^2 \int_\sigma f^2 ds.
\end{equation}
\end{conj}
It is straightforward to verify that this inequality holds when $\sigma$ is the round circle; in this case one has equality for the function $f=1$. 
\begin{rem}
This conjecture is termed the ``Oval's problem" and seems to have first appeared in the literature in \cite{BenLoss}.  In that paper, R. D. Benguria and M. Loss show that Conjecture \ref{OvalConj}
is related to conjectures about the one-dimensional Lieb-Thirring inequality.  They also 
prove that \eqref{BestConstIneq} holds if one replaces $(2\pi)^2$ by $ \frac{1}{2} (2\pi)^2$.  One of the 
difficulties in proving this conjecture seems to be that there is a whole family of curves on which the putative best 
constant $(2\pi)^2$ is achieved.  These were constructed by A. Burchard and L. E. Thomas  in \cite{Burchard} and consist 
of a one parameter family of ovals that contain the round circle and degenerate into a multiplicity 
two line segment.  Burchard and Thomas also show that in some neighborhood of the family they construct the conjecture holds.  For the general conjecture, the best constant so far achieved is $\approx 0.6 (2\pi)^2$ in \cite{Linde}.
\end{rem}

Using Conjecture \ref{OvalConj}, we show the following proposition which is a sharpening of Corollary \ref{OSCor}. The proof completely avoids the use of the Weierstrass representation.
\begin{prop}\label{AltProp}
Fix $\Omega$ the open region between two parallel planes $P_1=\set{x_3=h_1}$ and $P_2=\set{x_3=h_2}$ in $\Real^3$ where $h_1<h_2$.  Let $\Sigma\in \mathcal{A}(\Omega)$ and for $t\in (h_1, h_2)$ set $\Sigma_t=\Sigma \cap \set{x_3=t}$.
Then
\begin{equation*}
\frac{d^2}{dt^2} \mathcal{H}^1 (\Sigma_t)\geq \frac{(2\pi)^2}{F_3^2} \mathcal{H}^1 (\Sigma_t)+\int_{\Sigma_t} \frac{|A(\nu, E_2)|^2}{|\nabla x_3|^2} d\mathcal{H}^1.
\end{equation*}
Here $\nu, E_2$ are a global orthonormal frame on $\Sigma$ so that $\nu$ is parallel to $\nabla_\Sigma x_3$. 
\end{prop}
\begin{rem}
The term $A(\nu, E_2)$ at a point $p\in \Sigma$ measures the rate of change at $p$ of the ``contact angle" between $\Sigma$ and the plane $P=\set{x_3=x_3(p)}$ along the curve $\Sigma\cap P$.  Recall the contact angle at $p$ is the angle between $\mathbf{n}(p)$, the normal to $\Sigma$ at $p$, and the plane $P$.  In particular, this term vanishes identically on a vertical catenoid.  Indeed, the everywhere vanishing of such a term characterizes the vertical catenoid -- see \cite{Pyo}.
\end{rem}  

Before proving the proposition we do a slightly more general computation:
\begin{lem}\label{AltLem}
Consider $\Sigma$ a minimal hypersurface in $\Real^{n+1}$. Suppose that $\set{x_{n+1}=t}$ meets $\Sigma$ transversely for all $-\epsilon<t<\epsilon$ and that the intersection $\Sigma_t$ is a closed manifold.
Then
\begin{equation} \label{TwoDerivEqn}
\frac{d^2}{dt^2}\mathcal{H}^{n-1} (\Sigma_t)=\int_{\Sigma_t} \bigg|\nabla_{\Sigma_t}  \frac{1}{|\nabla_\Sigma x_{n+1}|}\bigg|^2 +\frac{(H_{\Sigma_t})^2+\sum_i |\beta_i|^2 +(H^{\Sigma}_{\Sigma_t})^2-|A^\Sigma_{\Sigma_t}|^2}{|\nabla_\Sigma x_{n+1}|^2} .
\end{equation}
Here $H_{\Sigma_t}$ is the mean curvature of $\Sigma_t$ as a codimension two surface in $\Real^{n+1}$ and $H^{\Sigma}_{\Sigma_t}$ is the mean curvature of $\Sigma_t$ as a hypersurface in $\Sigma$.  Similarly, $A^\Sigma_{\Sigma_t}$ is the second fundamental form of $\Sigma_t$ as a hypersurface in $\Sigma$.  Finally, 
\begin{equation*}
\beta_i=A(\nu, E_i)
\end{equation*} 
where $A$ is the second fundamental form of $\Sigma$, $\nu$ is a vector field on $\Sigma$ so $E_2, \ldots E_{n}$ are an orthormal frame on $\Sigma_t$ and $\nu$ is normal in $\Sigma$ to $\Sigma_t$. 
\end{lem}
\begin{proof}
The lemma will follow from the second variation formula for area. For $t\in (-\epsilon, \epsilon)$, let $\phi_t: \Real^{n+1}\to \Real^{n+1}$ denote a smooth family of $C^1$ diffeomorphisms of $\Real^{n+1}$ with $\phi_0(\xX)=\xX$ and $\phi_t$ equal to the identity outside of a compact set. Then we may write
\begin{equation*}
\phi_t(\xX)=\xX+t \mathbf{X}(\xX) +\frac{1}{2} t^2 \mathbf{Z}(\xX)+O(t^3)
\end{equation*}
where $\mathbf{X},\mathbf{Z}$ are compactly supported vector fields.   Fixing $M\subset \Real^{n+1}$ a $k$-dimensional compact surface and letting $M_t=\phi_t(M)$ the second variation formula (see \cite{SimonBook}) gives:
\begin{equation} \label{SecondVarEqn}
\left.\frac{d^2}{dt^2}\right|_{t=0} \mathcal{H}^k (M_t)=\int_{M} \Div_M \mathbf{Z}+(\Div_M \mathbf{X})^2+\sum_{i=1}^k| (D_{\tau_i} \mathbf{X})^\perp|^2-\sum_{i,j=1}^k \left( \tau_i \cdot D_{\tau_j} \mathbf{X}\right) \left( \tau_j \cdot D_{\tau_i} \mathbf{X}\right).
\end{equation}

We claim the lemma is a simple consequence of this formula.  Indeed,  for fixed $t_0\in (-\epsilon, \epsilon)$ let
$\mathbf{X}, \mathbf{Z}$ be vector fields normal to $\Sigma_{t_0}$ given by
\begin{equation*}
\mathbf{X}=\frac{\nabla_\Sigma x_{n+1}}{|\nabla_\Sigma x_{n+1}|^2}=\frac{1}{|\nabla_\Sigma x_{n+1} |} \nu
\end{equation*}
and
\begin{align*}
\mathbf{Z}&=-\frac{A(\nu, \nu) }{|\nabla_\Sigma x_{n+1}|^2}  \left(\mathbf{n} -\frac{\mathbf{n}\cdot \mathbf{e}_{n+1}}{|\nabla_\Sigma x_{n+1}|} \nu\right) \\
&=\frac{-H_\Sigma-\mathbf{H}_{\Sigma_t}\cdot \mathbf{n}}{|\nabla_\Sigma x_{n+1}|^2} \left(\mathbf{n} -\frac{\mathbf{n}\cdot \mathbf{e}_{n+1}}{|\nabla_\Sigma x_{n+1}|} \nu\right)\\
&=-\frac{H_\Sigma+\mathbf{H}_{\Sigma_t}\cdot \mathbf{n}}{|\nabla_\Sigma x_{n+1}|^3}\mathbf{N}.
\end{align*}
Here $\mathbf{n}$ is the normal to $\Sigma$, $\nu$ is the conormal to $\Sigma_{t_0}$ in $\Sigma$ and $\mathbf{N}$ is the outward normal to $\Sigma_{t_0}$ as a hypersurface in $\set{x_{n+1}=t_0}$.
Using these vector fields, given a parameterization $\FF_0$ of $\Sigma_{t_0}$ if we set
\begin{equation*}
\FF(\cdot, t)=\FF_0(\cdot)+ (t-t_0) \mathbf{X}(\cdot)+\frac{1}{2} (t-t_0)^2 \mathbf{Z}(\cdot)+\mathbf{G}(\cdot, t)
\end{equation*}
then $\FF(\cdot, t)$ is a parameterization of $\Sigma_t$ for $t$ near $t_0$ with
\begin{equation*}
\mathbf{G}(\cdot, t)=O(|t-t_0|^3).
\end{equation*}
In particular, \eqref{TwoDerivEqn} will follow from \eqref{SecondVarEqn} by using these vector fields.

It remains to evaluate the various terms in \eqref{SecondVarEqn}.  
We first compute:
\begin{align*}
\Div_{\Sigma_{t_0}} \mathbf{Z}&= -\mathbf{H}_{\Sigma_t} \cdot \mathbf{Z} \\
&=- \mathbf{H}_{\Sigma_t} \cdot \left(-\frac{H_\Sigma+\mathbf{H}_{\Sigma_t}\cdot \mathbf{n}}{|\nabla_\Sigma x_{n+1}|^3}\mathbf{N}\right)\\
&=\frac{ H_{\Sigma_t}^2}{|\nabla_{\Sigma} x_{n+1} |^2}
\end{align*}
where the last equality follows from the minimality of $\Sigma$.
One also computes:
\begin{align*}
\sum_{i=2}^{n} |(D_{E_i} \mathbf{X})^\perp|^2 &= \left|\nabla_{\Sigma_{t_0}} \frac{1}{|\nabla_\Sigma x_{n+1} |}\right|^2+\sum_{i=2}^n \frac{|A(E_i, \nu_{t_0})|^2}{|\nabla_\Sigma x_{n+1}|^2},
\end{align*}
\begin{equation*}
\sum_{i,j=2}^{n} \left(E_i \cdot D_{E_j} \mathbf{X}\right) \left( E_j \cdot D_{E_i} \mathbf{X} \right) =\frac{|A_{\Sigma_t}^\Sigma|^2}{|\nabla_\Sigma x_{n+1}|^2},
\end{equation*}
and 
\begin{equation*}
(\Div_{\Sigma_{t_0}} \mathbf{X})^2= \left( H^{\Sigma}_{\Sigma_t} \frac{1}{|\nabla_\Sigma x_{n+1}|}\right)^2.
\end{equation*}
Substituting these into \eqref{SecondVarEqn} completes the proof.
\end{proof}

We now show how Proposition \ref{AltProp} follows from Conjecture \ref{OvalConj}:
\begin{proof}
Set $\Sigma_t=\Sigma \cap \set{x_{3}=t}$.  By Lemma \ref{TransverseLem} all the $\Sigma_t$ are smooth curves.  As $\Sigma_t$ is a curve, $H_{\Sigma_t}=\kappa_{\Sigma_t}$, the geodesic curvature, and $|A_{\Sigma_t}^{\Sigma}|^2=(H^{\Sigma}_{\Sigma_t})^2$.  Thus, Lemma \ref{AltLem} gives:
\begin{align*}
\frac{d^2}{dt^2}\mathcal{H}^1 (\Sigma_t)&=\int_{\Sigma_t} \bigg|\nabla_{\Sigma_t}  \frac{1}{|\nabla_\Sigma x_{3}|}\bigg|^2 +\frac{\kappa_{\Sigma_t}^2+ |\beta_2|^2}{|\nabla_\Sigma x_{3}|^2} 
&\geq \frac{(2\pi)^2}{\mathcal{H}^1(\Sigma_t)^2} \int_{\Sigma_t}  \frac{1}{|\nabla_\Sigma x_{3}|^2}+  \int_{\Sigma_t}\frac{|\beta_2|^2}{|\nabla_\Sigma x_{3}|^2} .
\end{align*}
Here the inequality used Conjecture \ref{OvalConj}.
Set 
\begin{align*}
\alpha_t &=\frac{1}{\mathcal{H}^1(\Sigma_t)}\int_{\Sigma_t} \frac{1}{|\nabla_\Sigma x_3|} \geq \frac{1}{\mathcal{H}^1(\Sigma_t)}  \frac{ (\mathcal{H}^1(\Sigma_t))^2}{\int_{\Sigma_t} |\nabla_\Sigma x_3|}  = \frac{\mathcal{H}^1( \Sigma_t)}{F_3},
\end{align*}
where the second inequality follows from the Cauchy-Schwarz inequality and the last equality uses the fact that $\mathbf{e}_3\cdot \nu=\mathbf{e}_3\cdot \frac{\nabla_\Sigma x_3}{|\nabla_\Sigma x_3|}=|\nabla_\Sigma x_3|$.  Note that one has equality if and only if $|\nabla_\Sigma x_3|$ is constant on $\Sigma_t$.
Then on $\Sigma_t$:
\begin{equation*}
\frac{1}{|\nabla_\Sigma x_3|}=\alpha_t+\psi
\end{equation*}
where $\psi$ is a smooth function on $\Sigma_t$ with $\int_{\Sigma_t} \psi=0$.
Then one has:
\begin{align*}
\int_{\Sigma_t} \frac{1}{|\nabla_\Sigma x_3|^2} &=\alpha_t^2 \mathcal{H}^1 (\Sigma_t) +\int_{\Sigma_t} \psi^2 \geq \frac{(\mathcal{H}^1(\Sigma_t))^3}{F_3^2} +\int_{\Sigma_t} \psi^2.
\end{align*}
 Hence:
 \begin{equation*}
 \frac{d^2}{dt^2}\mathcal{H}^1 (\Sigma_t)\geq \frac{(2\pi)^2}{F_3^2}\mathcal{H}^1 (\Sigma_t)+  \int_{\Sigma_t}\frac{|\beta_2|^2}{|\nabla_\Sigma x_{3}|^2} .
 \end{equation*}
 \end{proof}
\bibliographystyle{amsplain}
\bibliography{Biblio}

\end{document}

%% file: MSCat2.tex

\begingroup
  \makeatletter
  \providecommand\color[2][]{%
    \errmessage{(Inkscape) Color is used for the text in Inkscape, but the package 'color.sty' is not loaded}
    \renewcommand\color[2][]{}%
  }
  \providecommand\transparent[1]{%
    \errmessage{(Inkscape) Transparency is used (non-zero) for the text in Inkscape, but the package 'transparent.sty' is not loaded}
    \renewcommand\transparent[1]{}%
  }
   \providecommand\rotatebox[2]{#2}
  \ifx\svgwidth\undefined
    \setlength{\unitlength}{384.8pt}
  \else
    \setlength{\unitlength}{\svgwidth}
  \fi
  \global\let\svgwidth\undefined
  \makeatother
  \begin{picture}(1,0.62577963)%
    \put(0,0){\includegraphics[width=\unitlength]{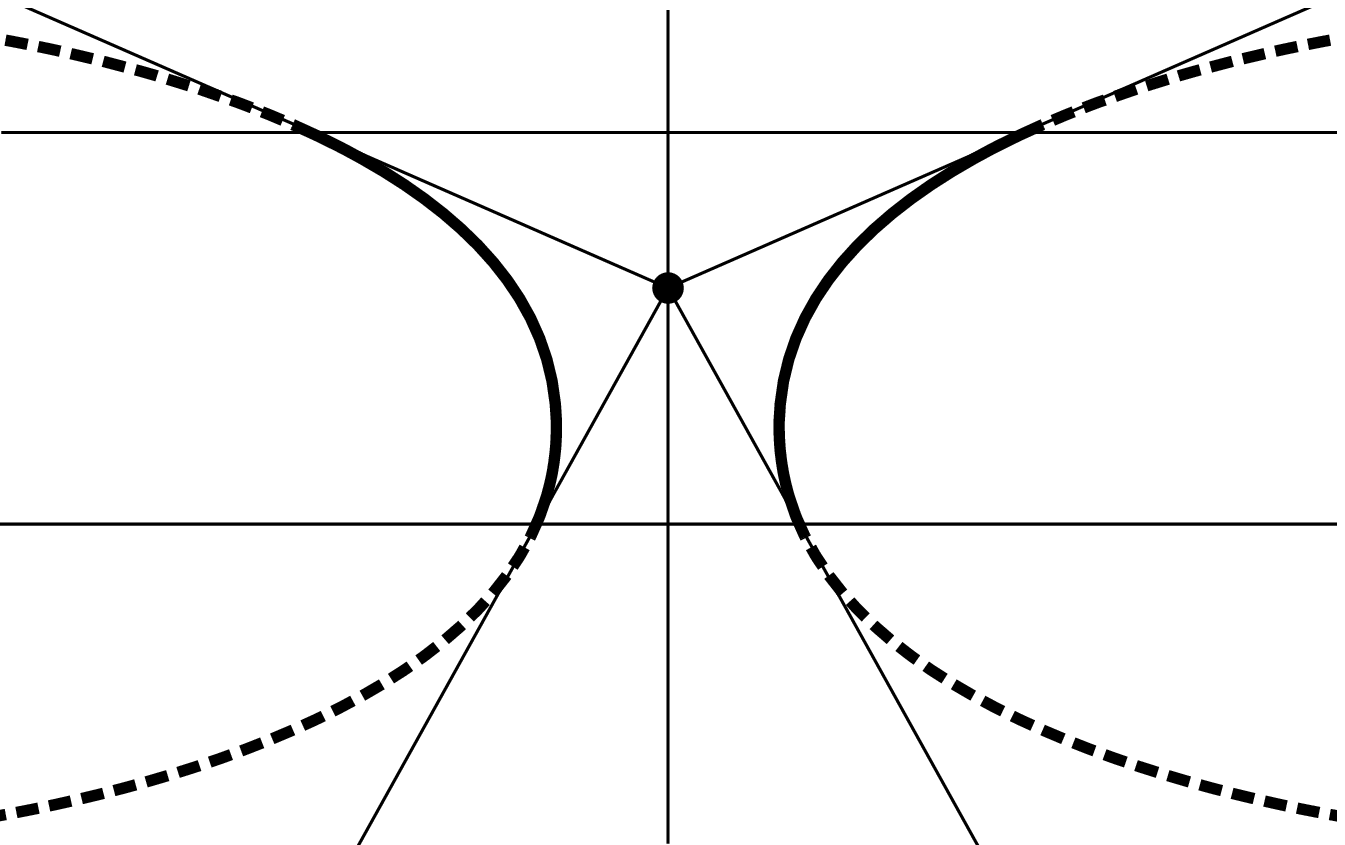}}%


\put(0.52076509,0.41){\makebox(0,0)[lt]{\begin{minipage}{0.05880306\unitlength}\raggedright $p$\end{minipage}}}%
\put(0.35554644,0.51){\makebox(0,0)[lt]{\begin{minipage}{0.15750823\unitlength}\raggedright $C_+(p)$\end{minipage}}}%
    \put(0.33365195,0.11739314){\makebox(0,0)[lt]{\begin{minipage}{0.15750823\unitlength}\raggedright $C_-(p)$\end{minipage}}}%
    \put(0.04260698,0.51239531){\makebox(0,0)[lt]{\begin{minipage}{0.20615861\unitlength}\raggedright $\set{x_3=t_+(p)}$\end{minipage}}}%
    \put(0.04260698,0.22133501){\makebox(0,0)[lt]{\begin{minipage}{0.20615861\unitlength}\raggedright $\set{x_3=t_-(p)}$\end{minipage}}}%
    \put(0.76506775,0.14458106){\makebox(0,0)[lt]{\begin{minipage}{0.24151262\unitlength}\raggedright $Cat$\end{minipage}}}%
    \put(0.60392923,0.38764428){\makebox(0,0)[lt]{\begin{minipage}{0.27301425\unitlength}\raggedright $Cat_{MS}(p)$\end{minipage}}}%
  \end{picture}%
\endgroup